\documentclass [12pt] {article}
\title {On homotopies with triple points of classical knots}
\author{Thomas Fiedler and Arnaud Mortier}
\usepackage{amssymb,amsfonts,epsfig,amsmath,graphics,graphicx}
\begin{document}
\newtheorem{proposition}{Proposition}
\newtheorem{theorem}{Theorem}
\newtheorem{lemma}{Lemma}
\newtheorem{corollary}{Corollary}
\newtheorem{example}{Example}
\newtheorem{remark}{Remark}
\newtheorem{definition}{Definition}
\newtheorem{question}{Question}
\newtheorem{conjecture}{Conjecture}
\maketitle
\begin{abstract}
We consider a knot homotopy as a cylinder in 4-space. An ordinary triple point $p$ of the cylinder is called {\em coherent}
if all three branches intersect at $p$ pairwise with the same index. A {\em triple unknotting} of a classical knot $K$ is a
homotopy which connects $K$ with the trivial knot and which has as singularities only coherent triple points.

We give a new formula for the first Vassiliev invariant $v_2(K)$ by using triple unknottings. As a corollary we obtain a 
very simple proof of the fact that passing a coherent triple point always changes the knot type. As another corollary we 
show that there are triple unknottings which are not homotopic as triple unknottings even if we allow more complicated 
singularities to appear in the homotopy of the homotopy.

\footnote{
classical knots, discriminant , triple unknottings}
\end{abstract}

\section{Introduction and results}
A generic homotopy between two classical knots in 3-space has only ordinary double points as singularities. Any two
knots can be connected by a generic homotopy. Moreover, this homotopy with fixed end points is unique up to homotopy,
because the space of all (including singular) knots is a contractible space.

But it turns out, that there are first degenerations of a generic homotopy which are already interesting. In this paper we
study triple unknottings as defined in the abstract. Every knot admits a triple unknotting as was proven in \cite{F1}. 

{\em But are triple unknottings also unique up to homotopy?}

In order to give a precise sense to this question we have to introduce the spaces with which we are working.

Let $\mathbb{R}^4$ be the euclidean space with coordinates $(x,y,z,t)$. Let $K_0$ be a non singular knot (i.e. an 
embedded circle) in $(x,y,z,0)$ and let $K_t,t \in [0,1]$ be a homotopy of $K_0$. We view the homotopy as a cylinder
in $(x,y,z,t)$. For each fixed $t$ we denote the orthogonal projection of $(x,y,z,t)$ to $(x,y,t)$ by $pr$. Hence, we can
consider a homotopy as a 1-parameter family of knot diagrams with respect to $pr$.

Let  $\Sigma$ be the discriminant of all {\em singular} knots in 3-space (i.e. the image $f(S^1) \subset \mathbb{R}^{3}$ is
not a submanifold). $\Sigma$ has a natural stratification (compare \cite{V}). The strata of codimension one 
$\Sigma^{(1)} \subset \Sigma$ are formed by all knots  which have exactly one ordinary double point. 
Let $\Sigma^{(2)}_{triple}$ be the union of all those strata of codimension two which correspond to knots 
which have exactly one ordinary triple point as singularity.

\begin{definition}
A homotopy of a knot $K_0$ is called a {\em triple homotopy} if it  intersects  $\Sigma$ only in 
$\Sigma^{(2)}_{triple}$,  where the intersection is transverse and coherent.
It is called {\em regular} if moreover the isotopy part of it is regular (i.e. involves no Reidemeister move of type I).

\end{definition}

\begin{remark}
 Notice that being coherent is a property of the  {\em homotopy} and not only of the stratum it goes across.
\end{remark}

\begin{definition}

A {\em homotopy of a triple homotopy} is a homotopy $h_s, s\subset [0,1]$ with fixed end points and such that $h_s$ is
 a triple homotopy for fixed $s$ besides 
for a finite number of $s$ where it either intersects $\Sigma$ in an adjacent stratum of codimension at least 3 or it becomes 
tangential to $\Sigma^{(2)}_{triple}$, or it intersects $\Sigma$ in the intersection of two strata of $\Sigma^{(2)}_{triple}$ (see Fig. 1).
\end{definition}

%1
\begin{figure}
\centering 
\psfig{file=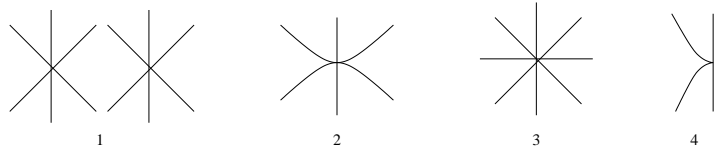}
\caption{the four strata we may go across when homotoping a homotopy}
\end{figure}

In other words, we add to the ordinary triple points all strata of higher codimension, but we still want that the homotopy
does not intersect $\Sigma^{(1)}$ and that it intersects $\Sigma^{(2)}_{triple}$ in a coherent way but perhaps not 
transversally.

It is easy to see that generically a coherent triple point in a homotopy corresponds to a 1-parameter family of diagrams 
as shown in Fig.2 and Fig.3 (compare \cite{F1}). We call the family in Fig.2 a {\em braid-like move} and those in Fig.3  a 
{\em star-like move}. We will see in the next section that braid-like moves and star-like moves can be replaced by each
other by a homotopy of a triple homotopy which passes through a triple point in an auto-tangency.

%2
\begin{figure}[h!]
\centering 
\psfig{file=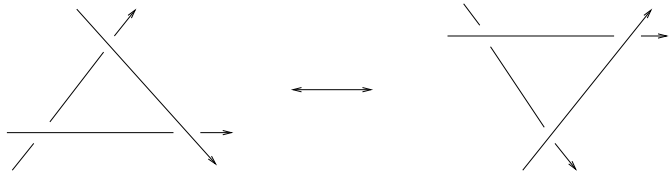}
\caption{braid-like moves}
\end{figure}

%3
\begin{figure}[h!]
\centering 
\psfig{file=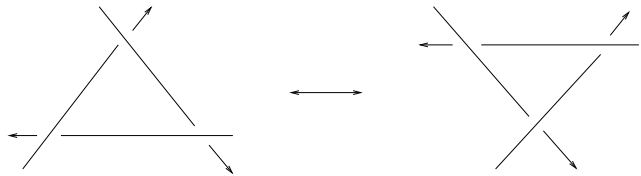}
\caption{star-like moves}
\end{figure}

The {\em writhe} of a diagram is defined as usual as the sum of the signs of the crossings (see e.g. \cite{K}). Notice that
the sum of the signs of the three involved crossings changes from -3 to +3 for both moves, braid-like or star-like.

Let $K_t$ be a triple homotopy and let $p$ be a coherent triple point in it.

\begin{definition}
 The index $ind(p)=+1$ if the writhe increases by +6 for increasing parameter $t$ and $ind(p)=-1$ otherwise.
The {\em index of the triple homotopy} $ind(K_t)$ is defined as the sum of $ind(p)$ over all triple points in $K_t$.
\end{definition}

\begin{theorem}
The index $ind(K_t)$ is invariant under homotopies of the triple homotopy $K_t$.
\end{theorem}

\begin{remark}
The index takes its values in an abelian group. It follows that $ind(K_t)$ is also invariant under homologies of triple 
homotopies (which are defined in the obvious way).
\end{remark}

We construct an example to answer the above question.

\begin{corollary}
There are two triple unknottings of the right trefoil which are not homotopic as triple unknottings.
\end{corollary}

\begin{corollary}
There are triple loops which are not homotopic -as triple homo\-topies- to regular triple loops.
\end{corollary}

As already mentioned above, we can assume that all triple points in our triple unknotting correspond to braid-like 
moves. Let $K$ be a knot with a triple point which corresponds to a braid-like move. We consider the {\em Gauss
diagram} associated to $K$ (compare \cite{PV} and also \cite{F2}). It contains a triangle corresponding to the triple
point $p$.

\begin{definition}
The {\em writhe of K with respect to p}, denoted by $W(p)$, is defined as the sum of the writhe of all arrows which
intersect the triangle in the Gauss diagram corresponding to p. Let $K_t$ be a triple homotopy. $W(K_t)$ is defined as 
the sum of 
$ind(p)W(p)$ over all triple points in $K_t$.
\end{definition}

Let $v_2(K)$ be the Vassiliev invariant of degre 2 (normalized to be 0 for the unknot and +1 for the trefoil) (see \cite{V}).

\begin{theorem}
For any triple homotopy $K_t$ which contains only braid-like coherent triple points,

$W(K_t)=v_2(K_1)-v_2(K_0)$. 

\end{theorem}

\begin{remark}
This formula is not true in the case of non coherent triple points.
The natural analog of this formula in the case of double points is the usual skein relation for $v_2$. 
\end{remark}

We will see in the next section that $W(p)$ is always odd for a braid-like move.

\begin{corollary}
Passing a coherent triple point always changes the knot type.
\end{corollary}

\begin{remark}
This corollary is not new. As far as we know, its first proof is contained in \cite{F1}. There it is shown that the Arf 
invariant changes always by passing a coherent triple point. But the proof uses Rokhlins 
$\mathbb{Z}/2\mathbb{Z}$-valued quadratic form associated
to an orientable characteristic surface in an orientable 4-manifold. The proof which uses Theorem 2 is much simpler.
\end{remark} 

The following question seems to us to be interesting.

\begin{question}
Can the next Vassiliev invariant $v_3(K)$ also be defined by using triple homotopies?
\end{question} 

 \section{Proofs}
Figure 4 shows how one can make a braid-like move into a star-like one, and conversely. 
The two diagrams with autotangencies correspond to the strata of codimension 3 that the homotopy goes across.
\newline

%4
\begin{figure}[h!]
\centering 
\psfig{file=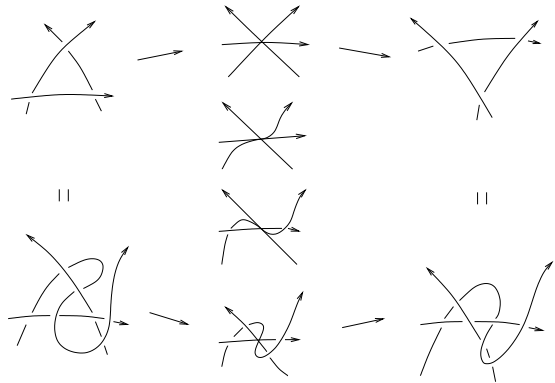}
\caption{a homotopy between braid and star moves}
\end{figure}

\begin{remark}
All the following arguments work as well with the second type of coherent braid-like triple points, obtained from 
figure 6 by reversing the orientation of any arrow involved in the middle branch (separating zone A and zone B).
\end{remark}

\textbf{Proof of Theorem 2.}
We use the expression of $v_2$ as a Gauss sum. Figure 5 reads as follows : to compute $v_2$ on a given knot $K$, 
draw its Gauss diagram, and choose a special point on it, not at the end of an arrow. 
Then for every couple of arrows whose configuration is as indicated by the picture, consider the product of the two writhes. 
The sum of all these products is $v_2(K)$, wherever the special point was chosen.

%5
\begin{figure}[ht!]
\centering 
\psfig{file=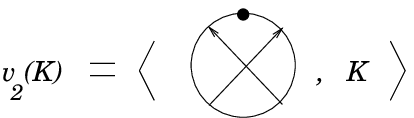}
\caption{$v_2$ as a Gauss sum}
\end{figure}

Now consider a triple homotopy $K_t$ which contains only one (braid-like) triple point (see Fig. 6). To compute $v_2(K_1)-v_2(K_0)$ 
we may restrict our attention to couples of arrows in which at least one arrow from the triple is implied. Let $w(X,Y)$ denote the sum 
of the writhes of all the arrows whose tail lies in zone $X$ and head in zone $Y$. The above description of $v_2$ yields :
\begin{align*} 
&v_2(K_1)-v_2(K_0)\\
&=(w(B,A)+w(B,C)+w(C,B))-(-w(A,C)-w(A,B)-w(C,A))\\
&= W(K_t). \\ 
\end{align*}
%6
\begin{figure}[h]
\centering 
\psfig{file=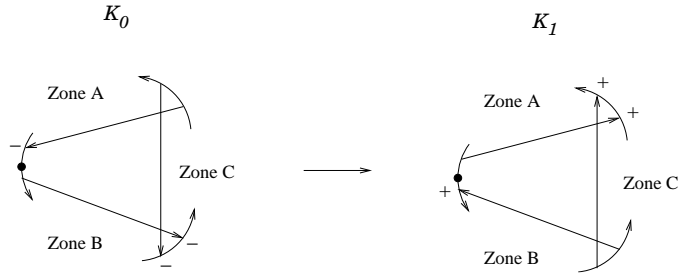}
\caption{an elementary triple homotopy}
\end{figure}

\textbf{Proof of Corollary 3.}
With the same notations as before, let us compute $W(K_t)$ again, this time the special point lying between the two 
crossings of the upper branch (which separates zone A and zone C). We get :
$$W(K_t) = 1+2w(B,A)+2w(B,C)+w(A,C)+w(C,A).$$
It remains to show that $w(A,C)+w(C,A)$ is even. Consider the knot $K_1$ with two crossings of the triple smoothed as 
indicated by figure 7. A, B and C have become the three components of a link, and we see that $w(A,C)$ and $w(C,A)$ are both 
equal to the linking number of A and B. It follows that $Ar\!f(K_0)\neq Ar\!f(K_1)$ -recall that $Ar\!f = v_2\,  mod\, 2$, which completes the proof.
\newline

%7
\begin{figure}[h]
\centering 
\psfig{file=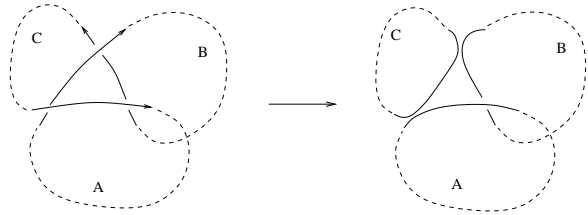}
\caption{smoothings}
\end{figure}

\textbf{Proof of Theorem 1.}
We have to prove that the index of each little loop around a stratum shown in figure 1 is zero.

Cases 1, 2 and 3 are easy to deal with : indeed, a generic little loop around one of these strata can be assumed to contain only
triple crossings and {\em regular} isotopies (which do not make use of $RI$). The writhe increase for such a path is 
linearly dependent on its index, which implies that a loop has index 0.
\newline

%8
\begin{figure}[h!]
\centering 
\psfig{file=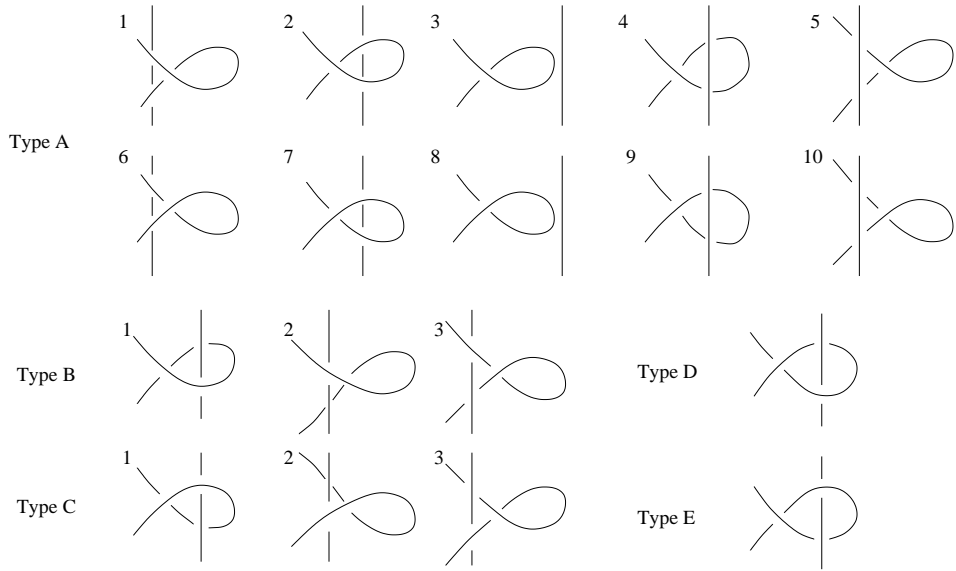}
\caption{resolutions of a cusped triple point}
\end{figure}

To treat case 4 (a cusped triple point), we show the 18 possible resolutions, divided into five different local knot types (see Fig. 8).
\begin{remark}
 We do not mention the orientations on Figure 8, because the picture does not depend on them.
\end{remark}

Diagrams of type A cannot contain a coherent triangle, therefore cannot be involved in a loop.
Figure 9 shows the {\em elementary} paths between the remaining diagrams, consisting either of going across a triple point, 
or of a little isotopy, for both possible orientations of the strands.

%9
\begin{figure}[h!]
\centering 
\psfig{file=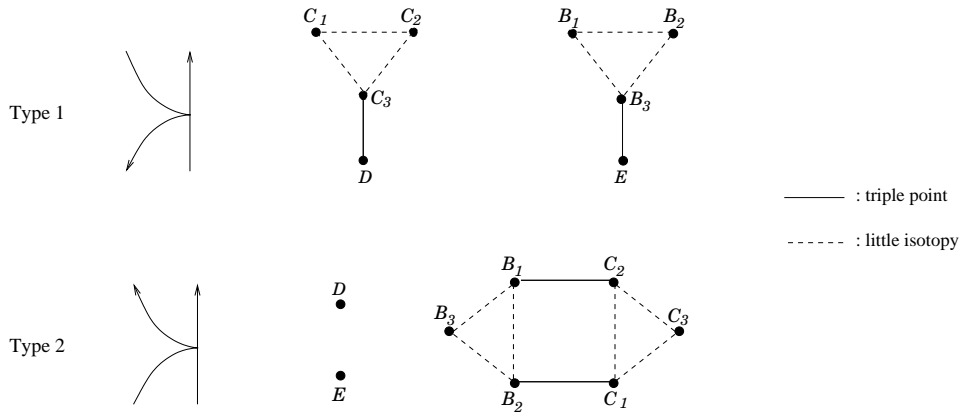}
\caption{the graph of little paths near a cusped triple point}
\end{figure}

The only loops appearing on the first graph of figure 9 consist of isotopies, and we may ignore them. The second graph 
contains one interesting loop, namely the square $B_1\rightarrow B_2\rightarrow C_1 \rightarrow C_2\rightarrow B_1$, 
whose index is clearly zero.

\begin{remark}
 Among the five knot types shown on figure 8, some might be equivalent (for instance in the case of a nugatory crossing), but
these unpredicted equivalences cannot be achieved by a little isotopy within the depicted part of $\mathbb{R}^3$.
\end{remark}

\textbf{Proof of Corollary 1.}
First, notice that if we know a triple homotopy between two knots $K_0$ and $K_1$, then for any knot $L$ we can deduce from 
it a triple homotopy from the connected sum $K_0\sharp L$ to $K_1\sharp L$ {\em with the same index}.

Figure 10 shows a positive elementary path from the trivial knot to the right trefoil $3_1^r$. On one hand, it proves that we 
may unknot $3_1^r$ with index -1. On the other hand, thanks to the connected sum property, it allows :
$$3_1^r \stackrel{+2}{\longrightarrow} 3_1^r\sharp 3_1^r\sharp 3_1^r.$$
Figure 11 shows a negative path from $3_1^r\sharp 3_1^r$ to $4_1$, from which we deduce :
$$3_1^r\sharp 3_1^r\sharp 3_1^r \stackrel{-1}{\longrightarrow} 4_1\sharp 3_1^r.$$
Figures 12 and 13 successively give :
$$4_1\sharp 3_1^r \stackrel{-1}{\longrightarrow} 3_1^l \stackrel{+1}{\longrightarrow} unknot.$$
Finally, we have two triple unknottings of $3_1^r$ with indices +1 and -1. But the index is a homotopy invariant by Theorem 1.
\newline

\textbf{Proof of Corollary 2.}
Corollary 1 implies the existence of a triple loop with index 2. But as we mentioned in the proof of Theorem 1, a regular 
loop must have index zero.

%10
\begin{figure}[ht!]
\centering 
\psfig{file=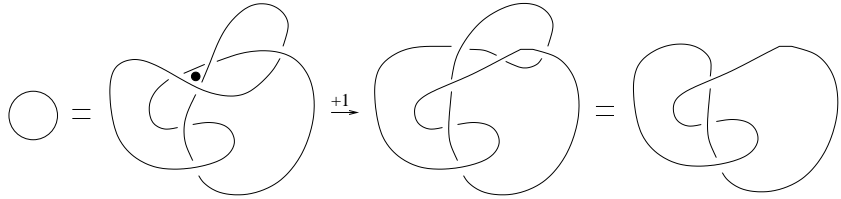}
\caption{trivial to trefoil}
\end{figure}

%11
\begin{figure}[ht!]
\centering 
\psfig{file=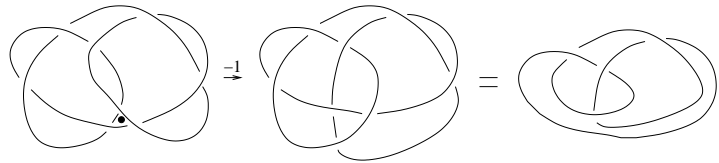}
\caption{trefoils to eight}
\end{figure}

%12
\begin{figure}[ht!]
\centering 
\psfig{file=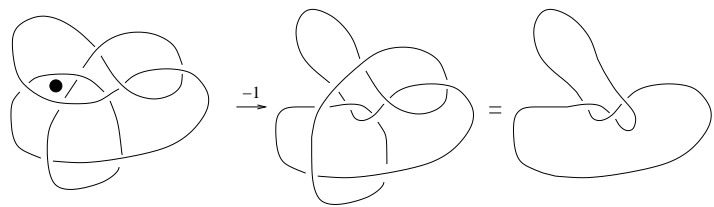}
\caption{eightfoil to trefoil}
\end{figure}

%13
\begin{figure}[ht!]
\centering 
\psfig{file=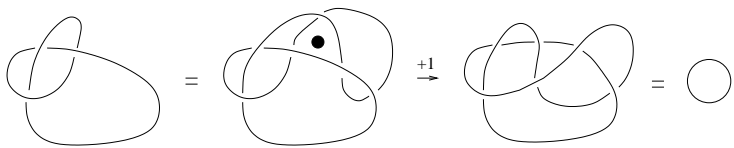}
\caption{left trefoil to trivial}
\end{figure}

Institut de Math\'ematiques de Toulouse

Universit\'e Paul Sabatier et CNRS (UMR 5219)

118, route de Narbonne 

31062 Toulouse Cedex 09, France

fiedler@picard.ups-tlse.fr

mortier@math.ups-tlse.fr

\end{document}